\documentclass[english, 12pt]{article}
\usepackage{amsfonts}
\usepackage{amsmath}
\usepackage{amssymb}
\usepackage{amsthm}
\usepackage{amscd}
\usepackage{amscd}
\usepackage{amsthm}

\usepackage{amsfonts}
\usepackage{amsmath}
\usepackage{color}
\usepackage{mathrsfs}

\newtheorem{defi}{Definition}[section]
\newtheorem{theo}{Theorem}[section]
\newtheorem{prop}{Proposition}[section]
\newtheorem{lem}{Lemma}[section]

\newtheorem{rk}{Remark}[section]

\newtheorem{exa}{Example}[section]
\newtheorem{coro}{Corollary}[section]
\numberwithin{equation}{section}

\def\calE{{\cal{E}}}

\def\calE{\mathcal{E}}

\def\R{\mathbb{R}}
\def\N{\mathbb{N}}
\def\Z{\mathbb{Z}}
\def\alp{\alpha}
\def\Om{\Omega}
\def\lam{\lambda}
\def\E{\mathbb{E}}
\def\P{{\bf P}}
\def\kmu{K^\mu}

\newcommand{\calD}{{\cal{D}}}

\DeclareMathOperator{\ran}{ran}
\DeclareMathOperator{\dom}{dom}

\title{Continuity aspects for traces of Dirichlet forms}
\usepackage{fancyhdr}
\makeatletter
\let\runtitle\@title
\makeatother
\begin{document}
\title{Continuity aspects for traces of  Dirichlet forms with respect to  monotone weak convergence of G-Kato measures}
\author{\normalsize
Ali BenAmor\footnote{High school for transport and logistics, University of Sousse, Sousse, Tunisia. E-mail: $ali_{-} ben05@yahoo.de$} 
}

\date{}

\maketitle

\begin{abstract}
We investigate some analytic properties of  traces of Dirichlet forms with respect to measures satisfying Hardy-type inequality. Among other results we prove convergence of spectra, ordered eigenvalues, eigenfunctions as well as convergence of resolvents on appropriate spaces, for traces of Dirichlet forms when the speed measure is the monotone weak limit of G-Kato measures. Some quantitative estimates are also given. As applications we show continuity of stationary solutions of some elliptic operators with measure-valued coefficients with respect to the coefficients and give an approximation procedure for the eigenvalues of one-dimensional graph Laplacian as well as for the Laplacian on annular thin sets with mixed Neumann-Wentzell boundary conditions.
\end{abstract}

\section{Introduction}

Given a regular Dirichlet form $\calE$ with domain in some Lebesgue space $L^2(X,m)$ and another measure $\mu$, let $\calE^\mu$ be the trace of $\calE$ with respect to $\mu$. We will give a new characterization of $\calE^\mu$ and establish an explicit formula for its resolvent, provided the measure $\mu$ satisfies a Hardy-type inequality. The obtained results will then be exploited to study continuous dependence of spectral objects with respect to variations of the underlying measures. Actually, under some conditions which will be precised later, we shall discuss continuity aspects of the map $\mu\mapsto\calE^\mu$. Precisely, if $(\mu_n)_{n\in\N\cup\{\infty\}}$ is a family of finite  G-Kato measures  such that $\mu_n\uparrow \mu_\infty$ or  $\mu_n\downarrow \mu_\infty$ weakly, we will show that spectral objects (such as eigenvalues, Riesz eigenprojections, eigenvectors) related to $\calE^{\mu_\infty}$ can be approximated by the corresponding spectral objects of the sequence $(\mu_n)_n$. Furthermore, we will show that for large integers  $n$, the dimensions of the spectral projections are locally constant, obtaining thereby convergence of the ordered eigenvalues and hence of Rayleigh--Ritz quotients. We will also establish convergence of the related resolvents and potential operators on the space of Borel bounded functions. For all considered approximations we will give errors estimates which depend mainly on the uniform norm of the potential of the difference $|\mu_\infty - \mu_n|$.\\
\indent Let us stress that the trace forms  $\calE^{\mu_n}$ have domains in a sequence of  Lebesgue spaces, namely $L^2(X,\mu_n)$, so that the above-mentioned  convergence problems become delicate to handle. Though there  is a method for discussing the problem by means of Mosco convergence on varying Hilbert spaces (see \cite{Kuwae}), however it  does not lead to any quantitative information. To overcome this drawback, thanks to appropriate assumptions which are fulfilled in many concrete examples (see Section \ref{Applications}), we will transform  the considered problem to a problem on a fixed Banach space, namely the space of Borel bounded functions on $X$. With this strategy one can use classical results regarding convergence, approximation and error estimates for sequences of operators, spectra and eigenfunctions  on a fixed Banach space.\\ 
\indent Our strategy  relies mainly on two facts: first we prove compactness of the related resolvents on Lebesgue spaces and on the space of Borel bounded functions,  then we  show that the respective spectra coincide  and establish a bijective correspondence between the respective  eigenfunctions. Second, we prove that weak monotone convergence of the considered sequence of measures leads to uniform convergence of the corresponding potentials and hence convergence of the related potential operators on the Banach space of Borel bounded functions.\\
\indent Though our assumptions on $\calE$ and   $(\mu_n)_{n\in\N\cup\{\infty\}}$ are somewhat  strong, they lead to satisfactory results. In fact, they recover and complement the results from \cite[Section 2.7]{Kuwae} as well as  the one-dimensional case studied in \cite{Ehnes-Hambly} and \cite{Freiberg-Minorics}.\\
\indent Concerning construction and analysis  of traces of forms, we mention that there are many methods at hands \cite{CF,HZ,F,AT,BBST,Post}.\\
We quote that traces of Dirichlet forms are also called time changed Dirichlet forms in the literature  \cite{CF}. In fact the  semigroup related to $\calE^\mu$ considered on the 'quasi-support' of the measure $\mu$ is obtained from the semigroup of the starting form $\calE$ by an appropriate time change.\\
\indent When $\calE$ corresponds to an elliptic partial differential operator on an open sufficiently regular set $\Omega$ in a  Riemannian manifold and $\mu$ is the surface measure of the boundary of $\Om$ the trace form $\calE^\mu$ corresponds to the Dirihlet-to-Neumann operator. Owing to the latter situation, traces of Dirichlet forms can be seen as quadratic forms corresponding to abstract Dirichlet-to-Neumann operators. 

\section{Traces of Dirichlet forms, revisited}
Let $X$ be a locally compact separable metric space, $m$ a positive Radon measure on Borel subset of $X$ with full support and $\calE$ a regular symmetric Dirichlet  form with domain  $\calD\subset L^2(X,m)$. We denote by $\|\cdot\|$ the $L^2(X,m)$ norm and set $\calE[u]:=\calE(u,u).$ Let $\mu$ be an other positive Radon measure on Borel subsets of $X$ charging no sets having zero capacity. We set $(\cdot,\cdot)_\mu$ resp. $\|\cdot\|_\mu$ the scalar product resp. the norm of $L^2(X,\mu)$. We designate by $\mathscr{B}_b(X)$ resp. $\mathscr{B}_b^+(X)$ the spaces of Borel bounded and positive Borel bounded functions on $X$ endowed with the topology of uniform convergence which we denote by $\|\cdot\|_\infty$.\\
For concepts related to Dirichlet forms we refer the reader to \cite{CF}.\\  
\indent We assume that $\calE$ is transient (this assumption is not restrictive. For, if $\calE$ were not transient one changes $\calE$ by $\calE_\alpha:=\calE + \alp\|\cdot\|^2$ for some  $\alpha>0$ which is always transient). Let  $\calD^e$  be the extended space of $\calE$ and $\calE^e$ the extension of $\calE$ on $\calD^e$. It is well known that $(\calD^e,\calE^e)$ is a Hilbert space. Moreover  elements from $\calD^e$ have $m$-representatives which are quasi continuous (q.c. for short) and two q.c. representatives coincide quasi everywhere (q.e. for short) and hence $\mu$ a.e., by assumption on $\mu$. Thus we will assume that elements from $\calD^e$ are q.c..\\
\indent We say that $\mu$ has finite energy integral whenever there is a constant $c>0$ such that
\[
\int_X |u|\,d\mu \leq c\sqrt{\calE^e[u]},\ u\in\calD^e.
\]
For a signed Radon measure $\nu$ charging not sets of zero capacity we say that  $\nu$ has finite energy integral if $\nu^\pm$ have finite energy integral.\\
\indent We assume that the following Hardy-type inequality holds: there is a constant $C_H>0$ such that
\begin{align}
\int_X u^2\,d\mu \leq C_H\calE^e[u],\ u\in\calD^e.
\label{hardy}
\end{align}
Hardy's inequality together with Riesz representation theorem imply that for every $u\in L^2(X,\mu)$ the signed measure $u\mu$ has finite energy integral and hence has  an $\alpha$-potential for any $\alpha\geq 0$ which we denote by $U_\alpha^\mu u$. For $\alpha = 0$  we omit  the subscript $0$. Moreover  $U_\alpha^\mu u \in\calD$ for any $\alpha>0$ and $U^\mu u\in\calD^e$.\\
For every $\alpha\geq 0$ and every $u\in L^2(X,\mu)$ the potential $U_\alpha^\mu u$ is characterized by being the unique element from $\calD$ ($\calD^e$ for $\alpha =0$) such that
\begin{eqnarray}
\calE_\alpha(U_\alpha^\mu u,v) = \int_X uv\,d\mu,\ \forall\,v\in\calD,\  \calE^e(U^\mu u,v) = \int_X uv\,d\mu,\  \forall\,v\in\calD^e.
\label{Potential}
\end{eqnarray}
We define the linear operators $K_\alpha^\mu, \alpha\geq 0$ by  
\[
\kmu_\alpha\colon\,L^2(X,\mu) \to L^2(X,\mu),\ \kmu_\alpha u = U_\alpha^\mu u.
\]
For $\alp=0$ we write $K^\mu$ instead of $K_0^\mu$.\\
\begin{rk}
{\rm 
By \cite[Theorem 3.1]{BA}, Hardy inequality holds  if and only if for all $u\in L^2(X,\mu)$ the signed measure $u\mu$ has finite energy integral and the operator $K^\mu$ is bounded on $L^2(X,\mu)$. Moreover, the best constant in the inequality coincides with $\|K^\mu\|$. 
}	
\label{Hardy-pot}
\end{rk}
Let us consider the linear operator
\[
J:\dom (J) =  (\calD^e,\calE^e)\to L^2(X,\mu),\ u\mapsto u.
\]
Owing to the fact that the measure $\mu$ does not charge sets of zero capacity, the map $J$ is well defined. Moreover, $J$ is bounded by Hardy's inequality  and has dense range by the regularity of $\calE$. Let us observe that
\begin{equation}
K^\mu = J U^\mu.
\label{potential-J}
\end{equation}
%
%
\begin{lem}
It holds $K^\mu = JJ^*$ and $K^\mu$ is an injective bounded selfadjoint positive operator.
\label{injective-potential}
\end{lem}
\begin{proof}
Let $u\in L^2(X,\mu), v\in\calD^e$. Then
\[
\calE^e(J^* u, v) = (u,Jv)_\mu = \int_X uv\,d\mu = \calE^e(U^\mu u,v).
\]	
Thus $J^* = U^\mu$. Using  formula (\ref{potential-J}) we obtain $JJ^* = JU^\mu = K^\mu$.\\
Let $u\in L^2(X,\mu)$ be such that $K^\mu u = 0$. By the first step we get $JJ^* u = 0$ and then
\[
\calE^e(J^* u, J^*u) = (u,JJ^*u)_\mu = 0.
\]	
Hence $J^* u = 0$. Thereby $u\in \ker J^* = (\ran J)^\perp = \{0\}$, by properties of the map $J$. Accordingly,  $K^\mu$ is injective. The rest follows from the first step of the proof.
\end{proof}
Let us briefly  recall the concept of the trace of the Dirichlet form $\calE$. Let $P^e$ be the $\calE^e$ orthogonal projection onto the $\calE^e$-orthogonal of $\ker(J)$ and $\calE^\mu$ be the quadratic form defined by 
\[
\dom(\calE^\mu) = \ran (J),\calE^\mu[Ju] = \calE^e[P^e u],\ u\in\ran (J). 
\]
It is known that $\calE^\mu$ is a Dirichlet form (\cite[Theorem 6.2.1]{FOT} or \cite[theorem 5.2]{BBST}) and is called the trace of $\calE$ w.r.t. the measure $\mu$.\\ Let $(R_\alpha)_{\alpha>0}$ be the resolvent family associated to $\calE^\mu$ and  $L_\mu$ be the positive self-adjoint operator related to $\calE^\mu$ via Kato representation theorem. By the standard theory relating resolvent families to positive self-adjoint operators, we have
\begin{eqnarray}
	\dom (L_\mu) = \ran (R_\alpha) = \ran (K^\mu),\ L_\mu u = R_\alpha^{-1}u - \alpha u = (K^\mu)^{-1} u,\ u\in  \ran (K^\mu).
	\label{Operator}
\end{eqnarray}
Furthermore
\begin{align}
\dom (\calE^\mu) = \ran \big((K^\mu)^{1/2}\big),\ \calE^\mu (u,v) = \big( (K^\mu)^{-\frac{1}{2}}u, (K^\mu)^{-\frac{1}{2}}v \big)_\mu,\ 
	u,v\in \ran \big((K^\mu)^{1/2}\big).
	\label{Form-formula}
\end{align}
\begin{prop}
The following assertions hold:
\begin{enumerate}
\item A resolvent formula:	
\begin{equation}
R_\alpha= (1 + \alpha\kmu)^{-1} \kmu,\ \alpha>0.
\label{resolvent-formula}
\end{equation}
\item $L_\mu$ has compact resolvent if and only if $K^\mu$ is compact.
\item All eigenvalues of $L^\mu$ are  strictly positive. Moreover $E$ is an eigenvalue of $L_\mu$ if and only  $1/E$ is an eigenvalue of $K^\mu$ with  the same eigenfunction and the same multiplicity.  
\end{enumerate}
\label{Resolvent-1}
\end{prop}
\begin{proof}
Assertion (2) is an immediate consequence of (1); while assertion (3) follows from (1) together with Hardy's inequality and positivity together with  injectivity of $K^\mu$.\\
(1): Let $\alpha>0$. By Lemma \ref{injective-potential} the operator $(1 + \alpha\kmu)^{-1} \kmu$ is everywhere defined and bounded on $L^2(X,\mu)$. Let $\psi\in L^2(X,\mu)$. We have to show that there is $u\in\calD^e$ such that  $Ju = (1 + \alpha\kmu)^{-1} \kmu \psi$ and 
\[
\calE^\mu(Ju,Jv) + \alpha\int_X uv\,d\mu = \int_X \psi v\,d\mu,\ \forall\,v\in\calD^e.
\]
Set $u:= U^\mu(1 + \alpha\kmu)^{-1}\psi$. Obviously $u\in\calD^e$ and
\[
Ju = JU^\mu(1 + \alpha\kmu)^{-1}\psi= K^\mu(1 + \alpha\kmu)^{-1}\psi = (1 + \alpha\kmu)^{-1}K^\mu\psi.
\]
On the other hand for every $v\in\ker(J)$ we have
\[
\calE^e(u,v) =   \calE^e(U^\mu(1 + \alpha\kmu)^{-1}\psi,v) = \int_X (1 + \alpha\kmu)^{-1}\psi\cdot v\,d\mu =0.
\]
Thus $P^e u = u$. Owing to the definition of $\calE^\mu$ we obtain, for every $v\in\calD^e$
\begin{align*}
\calE^\mu(Ju,Jv) &+ \alpha\int_X uv\,d\mu = \calE^e(P^e u,P^ev) + \alpha\int_X uv\,d\mu\\
& = \calE^e(P^e u,v) + \alpha\int_X uv\,d\mu=\calE^e(u,v) + \alpha\int_X uv\,d\mu\\
& =\int_X (1+\alpha K^\mu)^{-1}u\cdot v\,d\mu + \alpha\int_X K^\mu(1+\alpha K^\mu)^{-1}u\cdot v\,d\mu=\int_X uv\,d\mu,
\end{align*}
which completes the proof.
\end{proof}
For further investigations of the relationship between $\calE^e$ and $\calE^\mu$ let us give first a description of $L_\mu$ in term of the extended form $\calE^e$. To this end, for any function $u\in L^2(X,\mu)$ we denote by $\tilde u$ any function from $\calD^e$ such that $u= \tilde u\ \mu-a.e.$.
\begin{theo}
We have the characterization  of $L_\mu$:
\begin{align}
\dom (L_\mu) & = \big\{  u\in L^2(X,\mu)\colon\, \text{there is }\ \tilde u\in\calD^e\ \text{and}\  
 v\in  L^2(X,\mu)\nonumber\\ 
 & \text{s.t.}\ \calE^e (\tilde u, U^\mu w) = (v,U^\mu w)_\mu,\ \forall\,w\in  L^2(X,\mu) \big\},\nonumber\\
\label{domain}
& L_\mu u = v = (\kmu)^{-1} u.
\end{align} 
\label{OperatorDesc}
\end{theo}
\begin{proof}
We first claim that for any fixed $u,w\in L^2(X,\mu)$ the value of $\calE^e(\tilde u, U^\mu w)$ does not depend on the choice of $\tilde u$. Indeed,  if $\tilde u_1, \tilde u_2\in\calD^e$ are such that $\tilde u_1 - \tilde u_2 =0\ \mu-a.e.$ then
\[
\calE^e (\tilde u_1 - \tilde u_2 , U^\mu w) = \int_X  (\tilde u_1 - \tilde u_2) w\,d\mu = 0, 
\]
and the claim is proved.\\
We already know that $\dom (L_\mu) = \ran (\kmu)$. Let $u= \kmu v\in \dom (L_\mu)$. Then $v=L_\mu u$ and  we may and will choose $\tilde u = U^\mu v\in\calD^e$.  Let $w\in L^2(X,\mu)$. Then
\[
\calE^e(\tilde u, U^\mu w) = \calE^e( U^\mu v, U^\mu w) = \int_X v U^\mu w\,d\mu = (v,U^\mu w)_\mu.
\]
Conversely, let $V$ be the set on the right-hand-side of (\ref{domain}). We claim  that the operator given by the right-hand-side of  (\ref{domain}) is well defined. Indeed, let  $v_1,v_2\in L^2(X,\mu)$ be such that 
\[
\big(v_1 - v_2, U^\mu w \big)_\mu = 0,\ \forall\,w\in  L^2(X,\mu).
\]
Then
\[
\big(v_1 - v_2, U^\mu w \big)_\mu = 0 = \big(v_1 - v_2, K^\mu w \big)_\mu 
= \big(\kmu(v_1 - v_2),  w \big)_\mu = 0,\ \forall\,w\in  L^2(X,\mu).
\]
Thus $\kmu(v_1 - v_2) = 0$.  Since $\kmu$ is injective by Lemma \ref{injective-potential}, we obtain $v_1=v_2$.\\
Let $u\in V$. Then there is $\tilde u\in\calD^e$ and  $v\in L^2(X,\mu)$ such that 
\[
\calE^e(\tilde u, U^\mu w) = (v,U^\mu w)_\mu,\ \forall\, w\in L^2(X,\mu).
\]
Accordingly,
\begin{align*}
\calE^e(\tilde u, U^\mu w) &= (v,U^\mu w)_\mu = (v,K^\mu w)_\mu = (\kmu v, w)_\mu = (U^\mu v, w)_\mu\\
& = \calE^e(U^\mu v,U^\mu w)_\mu,\ \forall\, w\in L^2(X,\mu).
\end{align*}
Thus 
\[
\calE^e( \tilde u - U^\mu v, U^\mu w) = 0 = \big(\tilde u - U^\mu v, w \big),\ \forall\, w\in L^2(X,\mu), 
\]
leading to  $ \tilde u - U^\mu v =0\ \mu-a.e.$. Thereby $u= \kmu v$, $u\in \ran (\kmu) = \dom (L_\mu)$ and  $v= (\kmu)^{-1} u = L_\mu v$.
\end{proof}
We are now in position to characterize the form $\calE^\mu$. We set $\calD^\mu:=\dom(\calE^\mu)$.
\begin{theo}
\begin{enumerate}
\item The extended Dirichlet form $\calE^e$ with domain $\ran (U^\mu):=\{U^\mu u\colon\,u\in L^2(X,\mu)\}$ is closable in $L^2(X,\mu)$.
\item $\calE^\mu=\overline{\calE^e|_{\ran (U^\mu)}}$.
\end{enumerate}	
\label{Dirichlet-F}
\end{theo}
\begin{proof}
1. Let $(u_k)\subset L^2(X,\mu)$ be such that $\calE^e[U^\mu u_k -  U^\mu u_j]\to 0$ and  $ \|U^\mu u_k\|_\mu\to 0$. By the transience of $\calE$, the space $(\calD^e,\calE^e)$ is a Hilbert space. Thus there is $v\in \calD^e$ such that $\calE^e[U^\mu u_k - v]\to 0$. Using Hardy inequality we get $v=0\ \mu-a.e.$. Besides, an elementary computation leads to
\[
\calE^e[U^\mu u_k - v] = \calE^e[U^\mu u_k] -2\int_X U^\mu u_k v\,d\mu  + \calE^e[v]  \to 2\calE^e[v]=0.
\]
Thus $\calE^e[U^\mu u_k]\to 0$, which ends the proof of closability.\\
2. Since  $K^\mu u = U^\mu u\ \mu-a.e.$ for any $u\in L^2(X,\mu)$ and $\ran (K^\mu)$ is a core for $\calE^\mu$,  it suffices to show that $\calE^\mu[K^\mu u]=\calE^e[U^\mu u]$ for all $u\in L^2(X,\mu)$. Let $u\in L^2(X,\mu)$, then by (\ref{Form-formula})
\[
\calE^\mu[K^\mu u] = (K^\mu u,u)_\mu =  (U^\mu u,u)_\mu = \calE^e[U^\mu u].
\] 
\end{proof}
\section{Continuity aspects under monotone weak convergence of $G$-Kato measures}
Let $(\P_t)_{t>0}$ be the transition semigroup associated with $\calE$, i.e.,
\[
\P_t u(x) = \E^x(u(\mathbb{X}_t)),\ t>0, x\in X, u\in\mathscr{B}_b(X).
\]
From now on we assume that there is a jointly measurable symmetric function 
\[
p_t(\cdot,\cdot): (0,\infty)\times X\times X\to (0,\infty)
\]
such that 
\[
\P_t u(x) = \int_X p_t(x,y)u(y)\,dm(y),\  t>0, x\in X, u\in\mathscr{B}_b(X).
\]
For any $\alpha\geq 0$ we designate by $G_\alpha$ the  $\alpha$-order resolvent kernel of $\calE$:
\begin{equation}
G_\alpha (x,y):= \int_0^\infty e^{-\alpha  t} p_t(x,y)\,dt,\ x,y\in X.
\label{laplace-transform}
\end{equation}
Let $L$ be the positive self-adjoint operator associated with $\calE$. Then for any $\alpha>0, u\in L^2(X,m)$ we have
\[
(L + \alpha)^{-1} u = \int_X G_\alpha(\cdot,y)u(y)\,dm(y), m-a.e..
\]
For $\alpha=0$ we set $G:=G_0$. By transience of $\calE$ we have $G_\alpha\not\equiv\infty$, moreover $G_\alpha>0$ for any $\alpha\geq 0$.\\
We further assume that
\begin{equation}
\text{for every fixed $y\in X$ the function $G(\cdot,y)$ is lower semi-continuous  (l.s.c.) on $X$}
\label{lsc}
\end{equation}
We denote by $C_0(X)$ the space of real-valued continuous functions on $X$, vanishing at infinity.
\begin{defi}
Let $\mu$ be  a positive Radon measure on $X$. We say that $\mu$ is a $G$-Kato measure if 
\[
G^\mu1:=\int_X G(\cdot,y)\,d\mu(y) \in C_0(X).
\]
\end{defi}
We denote by 
 \[
 \mathscr{K}_b(X) := \{ \mu\ \text{is a finite}\ G -\text{Kato measure}\}.
 \]
We name $\mathscr{K}_b(X)$ the set of finite $G$-Kato measures.\\
In the Appendix we will give some conditions which ensure a measure $\mu$ to be a $G$-Kato measure.\\ 
\indent Our main goal, in this section,  is to establish continuity of spectral objects (essentially eigenvalues and eigenfunctions) of the trace form $\calE^\mu$ as well as the related resolvents under specific changes  of the measure $\mu$ in case $\calE^\mu$ has compact resolvent. Precisely, we will prove continuity of the mentioned objects if the measure $\mu\in\mathscr{K}_b(X)$ varies weakly in a monotone way. Besides, quantitative error estimates will be established.\\
Owing to Proposition \ref{Resolvent-1}, analysing continuity  of spectral  objects  of $\calE^\mu$ with respect to $\mu$ amounts to analysing continuity of those properties for $K^\mu$ with respect to $\mu$.\\
\indent To reach our objectives, we first prove that measures from $\mathscr{K}_b(X)$ satisfy Hardy's inequality and collect some properties related to the operators $K^\mu$  and $U^\mu$ on the space $\mathscr{B}_b(X)$.\\ 
For any  $\mu\in\mathscr{K}_b(X)$ and $u\in\mathscr{B}_b(X)\cup\mathscr{B}^+(X)$ we set
\[
G^\mu u(x) := \int_X G(x,y)u(y)\,d\mu(y),\ x\in X.
\]
\begin{lem}
Let $\mu\in\mathscr{K}_b(X)$. Then 
\[
\int_X u^2\,d\mu \leq \|G^\mu 1\|_\infty \calE^e[u],\ u\in\calD^e.
\]
\label{Kato-Hardy}
\end{lem}
\begin{rk}
{\rm
Similar inequality is established in \cite[Theorem 2.1]{Stollmann-Voigt}.
}
\end{rk}
\begin{proof}
It suffices to show the inequality for positive functions from $\calD^e$. Let  $\mu\in\mathscr{K}_b(X)$. Since Radon measures charging no sets of zero capacity are smooth in the sense of \cite[p.83]{FOT}, by \cite[Theorem 2.2.4, p.85]{FOT} there is an increasing generalized nest of closed sets such that $\mu_n:= 1_{F_n}\mu$ has finite energy integral and $\mu(X\setminus \cup_n F_n)=0$. Let $0\leq u\in C_c(X)\cap\calD$. It follows that the measure $\nu_n:= u\mu_n$ has finite energy integral. Thus   
\begin{equation}
\int_X u\,d\nu_n = \int_X u^2\,d\mu_n = \calE^e(U^{\nu_n}, u)\leq \sqrt{\calE^e[U^{\nu_n}]}\cdot \sqrt{\calE^e[u]}.
\label{Step1}
\end{equation}
An elementary computation leads to
\[
U^{\nu_n}  =  \int_X G(\cdot,y)\,d\nu_n(y) = \int_X G(\cdot,y)u(u)\,d\mu_n(y)\ q.e..
\]
Thereby we obtain
\begin{align}
\calE^e[U^{\nu_n}] &= \int_X U^{\nu_n}\,d\nu_n = \int_X \big(\int_X G(x,y)u(u)\,d\mu_n(y)\big) u(x)\,d\mu_n(x)\nonumber\\
&\leq \|G^\mu 1\|_\infty\int_X u^2\,d\mu_n.
\label{Step2}
\end{align}
Inequality (\ref{Step1}) in conjunction with inequality (\ref{Step2}) lead to 
\begin{equation}
	\int_{F_n} u^2\,d\mu \leq \|G^\mu 1\|_\infty\cdot \calE^e[u].
	\label{Step3}
\end{equation}
Hence
\begin{equation}
\int_X u^2\,d\mu = \int_{F_n} u^2\,d\mu	 + \int_{X\setminus F_n} u^2\,d\mu \leq \|G^\mu\ 1\|_\infty\cdot \calE^e[u]
+ \int_{X\setminus F_n} u^2\,d\mu,\ \forall\,n\in\N. 
\end{equation}	
Passing to the limit and using the fact that $\mu(X\setminus F_n)\to 0$ we get the desired inequality for $u\in C_c(X)\cap\calD$, positive. For positive $u\in\calD$ the proof follows from the regularity of $\calE$ and the fact that $\calE_1$ convergent sequence has a $q.e.$ convergent subsequence. The same arguments lead to the inequality for positive $u$ from $\calD^e$.
\end{proof}
\begin{rk}
{\rm
	\begin{enumerate}
		\item Every measure from $\mathscr{K}_b(X)$ has finite energy integral.
		\item In Lemma \ref{Kato-Hardy} one can drop the finiteness assumption on $\mu$ and still get Hardy's inequality.
	\end{enumerate}
	
	}
	
\end{rk}

By Lemma \ref{Kato-Hardy}, measures from $\mathscr{K}_b(X)$ satisfy Hardy's inequality and an elementary computation leads to 
\[
U^\mu u =  \int_X G(\cdot,y)u(y)\,d\mu(y)\ q.e.,\ u\in L^2(X,\mu),
\]
Moreover
\[
K^\mu u =  \int_X G(\cdot,y)u(y)\,d\mu(y)\ \mu-a.e.,\ u\in L^2(X,\mu).
\]
\begin{lem}
The linear operator
\[
G^\mu:=\mathscr{B}_b(X)\to  \mathscr{B}_b(X) ,\ u\mapsto \int_X G(\cdot,y)u(y)\,d\mu(y),
\]
is bounded and  $G^\mu(\mathscr{B}_b(X))\subset C_0(X)$.
\label{Kato-property}
\end{lem}
\begin{proof}
The boundedness of $G^\mu$ is obvious.\\
It suffices to prove $G^\mu(\mathscr{B}^+(X))\subset C_0(X)$. Let $u\in\mathscr{B}^+(X) $, by l.s.c. property for $G$, using Fatou lemma  we get that $G^\mu u$ is l.s.c.. On the other hand, by the fact that $\mu$ is G-Kato, the first part of the proof and the formula $G^\mu u = \|u\|_{\infty} G^\mu 1 - G^\mu(\|u\|_{\infty} - u)$ we conclude that $G^\mu u$ is upper semi-continuous. Thus $G^\mu$ is continuous and the inequality $0\leq  G^\mu u\leq \|u\|_\infty G^{\mu}1$ yields $G^\mu u\in C_0(X)$. 
\end{proof} 
\begin{lem}
The operators $K^\mu, G^\mu$ and  $K^\mu|_{L^\infty(X,\mu)}: L^\infty(X,\mu)\to L^\infty(X,\mu)$ are compact.
\label{compact}
\end{lem}
\begin{proof}
Compactness of $K^\mu$ follows from \cite[Theorem 2.5]{BH} together with \cite[Remark 2.1, 3]{BH}.\\
Let us prove compactness of $G^\mu$. Let $(u_k)_k\subset \mathscr{B}_b(X)$ be such that $0\leq u_k\leq 1$. Then by Lemma \ref{Kato-property}, $v_k:= G^\mu u_k$ satisfies 
\begin{align}
v_k\in C_0(X)\ \text{and}\ 0\leq\sup_k v_k(x)\leq G^\mu 1(x),\ x\in X.
\label{Unif-Estim}
\end{align}
Thus $(v_k)_k$ is a sequence of continuous potentials which are specific minorants of the continuous potential $G^\mu 1(x)$, i.e., $G^\mu 1 - v_k$ is a potential (in fact it is the potential of the measure $(1-u_k)\mu$). It follows from  \cite[Proposition 3.1]{Hansen-Polar}, that the sequence $(v_k)_k$ is equicontinuous. As it is uniformly bounded, Ascoli's theorem implies that there is a subsequence which we still denote by  $(v_k)_k$ and $v\in C(X)$ such that $v_k\to v$  locally uniformly on $X$. We claim that $v \in C_0(X)$. Indeed,  we have
\[
0\leq  v(x) = \lim_{k\to\infty} G^\mu u_k(x)\leq G^\mu 1(x),\ x\in X,
\]
which shows that $v\in C_0(X)$, by assumptions on the measure $\mu$.\\
It remains to show uniform convergence of $(v_k)$ towards $v$. Let $\epsilon>0$ and $A\subset X$ a compact set such that $\|1_{A^c} G^\mu 1\|_\infty<\epsilon$. Then from the above inequality we get
\[
\|1_{A^c} v\|_\infty<\epsilon\ \text{and}\ \|1_{A^c} v_k\|_\infty<\epsilon,\ \forall\,k.
\]
By local uniform convergence of $(v_k)$ towards $v$, we obtain $\|1_A(v_k - v)\|_\infty<\epsilon$ for large $k$. Hence, for large $k$ we get
\begin{align*}
\|v_k - v\|_\infty	 \leq \|1_A(v_k - v)\|_\infty + \|1_{A^c} v\|_\infty + \|1_{A^c} v_k\|_\infty<3\epsilon.
\end{align*}	
Therefore 
\[
\lim_{k\to\infty}\sup_{x\in X} \|G^\mu u_k(x) - v(x)\|_\infty = 0.
\]
Finally the proof of compactness of $K^\mu|_{L^\infty(X,\mu)}$ follows the same lines as before, so we omit it.
\end{proof}
\begin{rk}
{\rm
Let $\lam$ be a strictly positive eigenvalue of $G^\mu$. Then, owing to compactness and $L^2$-symmetry of $G^\mu$ we conclude (as it is the case for compact self-adjoint operators) that  
\begin{align*}
\ker(G^\mu-\lam) = \ker(G^\mu-\lam)^n,\ \forall\,n\geq 1.
\end{align*}
Thus the geometric and algebraic multiplicities of $\lam$ coincide. Consequently, from now on we will simply consider 'multiplicity' of $\lam$ to mean geometric or algebraic multiplicity. 
}
\label{general-remark-1}
\end{rk}
For any linear operator $A$ we denote by $\sigma(A)$ its spectrum and by $\rho(A)$ its resolvent.
\begin{theo}
The following hold true:
\begin{enumerate}
\item  $\sigma(K^\mu)\setminus\{0\}$ consists of at most a countable set of strictly positive eigenvalues with finite multiplicity, accumulating possibly at zero.
\item $\sigma(K^\mu) = \sigma(K^\mu|_{L^\infty(X,\mu)})$. Moreover every eigenvalue of $K^\mu$ is an eigenvalue of  $K^\mu|_{L^\infty(X,\mu)}$ with the same multiplicity.
\item  $\sigma(G^\mu)\setminus\{0\}$ consists of at most a countable set of positive eigenvalues with finite multiplicity, accumulating possibly at zero.
\item Every eigenvalue of $K^\mu$ is an eigenvalue of $G^\mu$ with the same multiplicity; conversely every strictly positive eigenvalue of $G^\mu$ is an eigenvalue of $K^\mu$ with the same multiplicity.
\item Every eigenvalue of $K^\mu$ has an eigenfunction which belongs to  $C_0(X)$.
\end{enumerate}
\label{spectral}
\end{theo}
\begin{proof}
The first and third assertions follow from the standard theory of compact operators together with injectivity and positivity of $K^\mu$ and the  $L^2(X,\mu)$-symmetry of $G^\mu$.\\	
The second assertion follows from the fact that $L^\infty(X,\mu)$ is dense in $L^2(X,\mu)$, the continuity of the inclusion map $L^\infty(X,\mu)\subset L^2(X,\mu)$, the  $L^2(X,\mu)$-symmetry of $K^\mu|_{L^\infty(X,\mu)}$ and the results from \cite[Sections 5.2-5.3]{BH}.\\
Let $\lambda$ be an eigenvalue of $K^\mu$. By assertion 2., $\lambda$ is an eigenvalue of $K^\mu|_{L^\infty(X,\mu)}$. Thus there is $u\in L^\infty(X,\mu)\setminus\{0\}$ such that $K^\mu u =\lam u\ \mu-a.e.$. Let $\tilde u\in \mathscr{B}_b(X)$ be such that $\tilde u = u\ \mu-a.e.$. Then
\[
G^\mu\tilde u = K^\mu u = \lambda u = \lam\tilde u\ \mu-a.e.
\]
Multiplying the latter identity by $G$ and integrating w.r.t. to $\mu$ we get
\[
G^\mu(G^\mu\tilde u) = \lam G^\mu\tilde u.
\]
Set $v:= G^\mu\tilde u$. Then $v\in C_0(X)\setminus\{0\}$ and is an eigenfunction of $G^\mu$ with $\lam$ being the corresponding eigenvalue. Following the same lines we show that if $u_1,u_2$ are $L^2$-linearly independent eigenfunctions of $\lam$ then $G^\mu\tilde u_1, G^\mu\tilde u_2$ are $\mathscr{B}_b(X)$-linearly independent eigenfunctions of $G^\mu$ corresponding to the eigenvalue $\lambda$.\\
Now let $\lam$ be a strictly positive eigenvalue of $G^\mu$ with corresponding eigenfunction $u$. Obviously $u\in L^2(X,\mu)$ and $K^\mu u = \lam u\ \mu-a.e.$. It suffices to show that $u\not\equiv 0$ in $L^2(X,\mu)$. For, if $u=0\ \mu-a.e.$, then $G^\mu u = 0$ everywhere which contradicts the fact that $u$ is an eigenfunction corresponding to a strictly positive eigenvalue. Finally it is easy to show that  $\mathscr{B}_b(X)$-linearly independent eigenfunctions corresponding to a strictly positive eigenvalue are also  $L^2(X,\mu)$-linearly independent, which finishes the proof of assertions 4-5. 
\end{proof}
The eigenvalues of $K^\mu$ will be arranged in a decreasing way, counted as many times as their respective multiplicities. 
\[
\lam_0\geq \lam_1\geq\cdots\lam_k\geq\lam_{k+1}\cdots(\to 0).
\] 
\begin{rk}
{\rm 
Since we are going to use Riesz projections we are lead to extend the  operator $G^\mu$ to complex valued functions. We denote by $\mathscr{B}_b^{\mathbb C}(X)$ the set of Borel measurable bounded complex-valued functions on $X$. Besides, we denote by $G_{\mathbb C}^\mu$ the operator $G^\mu$ acting on $\mathscr{B}_b^{\mathbb C}(X)$:  
\begin{align*}
	G_{\mathbb C}^\mu : \mathscr{B}_b^{\mathbb C}(X) \to \mathscr{B}_b^{\mathbb C}(X), u\mapsto  \int_X G(\cdot,y)u(y)\,d\mu(y).
\end{align*}
Then $G_{\mathbb C}^{\mu}$ is still compact and has the same spectrum as $\sigma(G^\mu)$. Moreover
\begin{equation}
	\dim\ker(G^\mu - \lam) = \dim\ker(G_{\mathbb C}^\mu - \lam)\ \text{for any}\ \lam\in\R.
	\label{kernel-dim-1} 
\end{equation} 
Let $\Gamma$ be a positively oriented circle enclosing $\lam$ but no other elements from the spectrum and $Q^\mu$ the Riesz projection of $\lam$ with respect to $G_{\mathbb C}^\mu$:
\[
Q^\mu:= -\frac{1}{2i\pi} \int_\Gamma (G_{\mathbb C}^\mu -z)^{-1}\,dz.
\]
Owing to compactness and $L^2$-symmetry of $G_{\mathbb C}^\mu$ once again, in conjunction with spectral properties of compact operators (see \cite[Theorem 2.21, p.53]{Dowson}), we obtain
\begin{align*}
	\ran (Q^\mu) = \ker (G_{\mathbb C}^\mu -\lam). 	
\end{align*}
Using (\ref{kernel-dim-1}) we achieve
\begin{align}
	\dim\ker(G^\mu - \lam) = \dim\ker(G_{\mathbb C}^\mu -\lam) = \dim\ran (Q^\mu).
	\label{kernel-dim-2}
\end{align}
}
\label{general-remark-2}
\end{rk}
We finish this section by a result relating dimensions of ranges of spectral projections to dimensions of ranges of Riesz projections; it will be used later.\\
\begin{lem}
Let $0<a<b$. Assume that $a,b$ are not eigenvalues of $\calE^\mu$. Let $\Gamma$ be a positively oriented circle connecting $1/a$ and $1/b$. Let $P{(a,b)}$ be the spectral projection of $(a,b)$ with respect to $L_\mu$ and
\[
Q^{\mu} = -\frac{1}{2i\pi} \int_\Gamma (G_{\mathbb C}^\mu - z)^{-1}\,dz.
\]
Then
\[
\dim \ran\big(P(a,b)\big) = \dim \ran(Q^\mu).
\]
\label{dimension-equal}
\end{lem}
\begin{proof}
Let $\Gamma'$ be the image of $\Gamma$ under the map $z\mapsto 1/z$, i.e., $\Gamma'$ is the positively oriented circle joining $a$ and $b$. 

Let $\lam_k', k\in\{1,\cdots,N\}$ be the set of distinct eigenvalues of $K^{\mu}$ lying inside $\Gamma$, $M:=\sum_{k=1}^N \dim \ker(K^\mu - \lam_k')$ and $(u_1',\cdots,u_M')$ an orthonormal basis of $\oplus_{k=1}^N \ker(K^\mu - \lam_k')$. Then by the already proved spectral properties for $K^\mu$, we obtain 
\[
P(a,b) = \sum_{k=1}^M (\cdot,u_k')_\mu u_k'.
\]
Thus $\dim \ran(P(a,b)) = M$.\\
Let us recall that the $\lam_k'$ are the sole eigenvalues of $G^\mu$ and hence of $G_{\mathbb C}^\mu$ lying inside $\Gamma$. Let $\Gamma_k$ be a positively oriented small circle around $\lam_k'$. Set
\[
Q_k^{\mu} := -\frac{1}{2i\pi} \int_{\Gamma_k} (G_{\mathbb C}^{\mu} - z)^{-1}\,dz.
\] 
Then  $\ran(Q_k^{\mu}) = \ker(G_{\mathbb C}^\mu - \lam_k')$ and by the holomorphic functional calculus we have
\[
Q^{\mu} = \oplus_{k=1}^N Q^\mu_k. 
\]
Thus $\dim \ran(Q^\mu) = \sum_{k=1}^N \dim\ran (Q_k^\mu)$.\\
On the other hand, using Theorem \ref{spectral} in conjunction with Remark \ref{general-remark-2} we obtain 
\[
\dim \ran(Q_k^{\mu}) = \dim\ker(G_{\mathbb C}^\mu - \lam_k') =  \dim\ker(G^\mu - \lam_k') = \dim\ker(K^\mu - \lam_k').
\]
Finally we obtain
\begin{align*}
\dim \ran (Q^\mu) &=  \sum_{k=1}^N \dim\ran (Q_k^\mu)= \sum_{k=1}^N \dim\ker(K^\mu - \lam_k') = M\\
& = \dim \big(\ran P(a,b)\big).
\end{align*}
which completes the proof.
\end{proof}

{\subsection{Qualitative results}}
We set $\N_\infty := \N\cup\{\infty\}$.\\
Let $(\mu_n)_{n\in\N_\infty}\subset \mathscr{K}_b(X)$. For each $n\in\N_\infty$, let $(\lam_n^{(k)})_{k\in I_n}$ be the eigenvalues of $K^{\mu_n}$, arranged as before. We would like to analyse convergence of eigenvalues, eigenvectors and resolvents of $\calE^{\mu_n}$ towards those of $\calE^{\mu_\infty}$ under appropriate convergence of measures, namely monotone weak convergence of measures. The discussion in the beginning of the section together with Theorem \ref{spectral} imply that investigation of continuous dependence of the mentioned objects related to $\calE^{\mu_n}$ with respect to $\mu_n$ amounts to investigate continuous dependence of those objects related to $G^{\mu_n}$ with respect to $\mu_n$ on the fixed Banach space $\mathscr{B}_b(X)$.\\
\indent Our assumption is: 
\begin{equation}
\mu_n\to  \mu_\infty\ \text{weakly-monotonically}.
\label{Weak-Monotone}
\end{equation}   
\indent{\em All the proofs in the remainder of the paper will be elaborated for $\mu_n\uparrow \mu_\infty$. Those corresponding to the decreasing case follow the same lines}.
\begin{theo}
Under assumption (\ref{Weak-Monotone}) the following hold:
\begin{enumerate}
\item If $\mu_n\uparrow \mu_\infty$ (resp. $\mu_n\downarrow \mu_\infty$ ) then $G^{\mu_n}1 \uparrow G^{\mu_\infty}1$ (resp.  $G^{\mu_n}1 \downarrow G^{\mu_\infty}1$ ), uniformly.
\item $\|G^{\mu_\infty} - G^{\mu_n}\| \leq \|G^{\mu_n}1 - G^{\mu_\infty}1\|_\infty$.
\end{enumerate}
\label{Uniform-Convergence}
\end{theo}
\begin{proof}
As $(G^{\mu_n} 1)_{n\in\N_\infty}\subset C_0(X)$, it suffices to show that $G^{\mu_n} 1 \uparrow G^{\mu_\infty} 1$ pointwise. By monotony of the $\mu_n$'s we have $0\leq G^{\mu_n} 1 \leq G^{\mu_\infty}1$ and then $\lim_{n\to\infty} G^{\mu_n} 1 (x) \leq   G^{\mu_\infty} 1 (x)$ for all $x\in X$.\\
Let $x\in X$ be fixed. Since $G(x,\cdot)$ is lower semi-continuous (by lower semi-continuity of $G(\cdot,y)$ and symmetry of $G$) and $G\not\equiv\infty$ there is a sequence $(G_n(x,\cdot))_n\subset C_b(X)$ such that  $G_n(x,y)\uparrow G(x,y)$ for all $y\in X$ (see \cite[Theorem 2.1.3, p.26]{Ransford}). Thus 
\[
G^{\mu_n} 1(x) = \int_X G(x,y)\,d\mu_n(y) \geq \int_X G_k(x,y)\,d\mu_n(y)\ \text{for all}\ k,x 
\]
and then by weak convergence of the $\mu_n$'s
\[
\lim_{n\to\infty}G^{\mu_n} 1(x) \geq \lim_{n\to\infty} \int_X G_k(x,y)\,d\mu_n(y)
 = \int_X G_k(x,y)\,d\mu_\infty(y)\ \text{for all}\ k,x. 
\]
Passing to the limit in the latter inequality and using monotone convergence theorem we get $\lim_{n\to\infty}G^{\mu_n} 1(x) \geq G^{\mu_\infty} 1(x)$ for all $x\in X$ and assertion 1. follows.\\
To prove the second assertion let $u\in\mathscr{B}_b(X)$ be such that $0\leq u\leq 1$.  Let us observe that $ G^{\mu_\infty} 1 - G^{\mu_n} 1 = G^{(\mu_\infty - \mu_n)}1$. Indeed, since $\mu_n\uparrow\mu_\infty$ there is a sequence of Borel bounded functions $(\varphi_n)$ such that $0\leq\varphi_n\leq 1$ and $\mu_n=\varphi_n\cdot\mu_\infty$, which leads to the claim, by an elementary computation. Thus
\begin{align*}
0\leq G^{\mu_\infty} u (x) - G^{\mu_n} u(x)\leq (G^{(\mu_\infty - \mu_n)} 1 )(x) \|u\|_\infty \leq \|G^{\mu_\infty} 1 - G^{\mu_n} 1\|_\infty.
\end{align*}	
Now let $u\in\mathscr{B}_b(X)$ be such that $\|u\|_\infty\leq 1$. We write $u=u^+ - u^-$ and use the latter inequality to get once again
\begin{align*}
|G^{\mu_\infty} u (x) - G^{\mu_n} u(x)|\leq \|G^{\mu_\infty} 1 - G^{\mu_n} 1\|_\infty,
\end{align*}
which leads to assertion 2. and the proof is completed.
\end{proof}

\begin{theo}
Under assumption (\ref{Weak-Monotone}) we have:
\begin{enumerate}
\item $\lim_{n\to\infty} \|G^{\mu_n} - G^{\mu_\infty}\| = 0$.
\item Let $\lam_\infty$ be a strictly positive eigenvalue of $G^{\mu_\infty}$ of multiplicity $m$. Then for large $n$ there is  $m$  eigenvalues (counted as many times as their respective multiplicities) $\Lambda_n^{(1)},\cdots\Lambda_n^{(m)}$ of $G^{\mu_n}$ such that  $\lim_{n\to\infty} \Lambda_n^{(j)} = \lam_\infty, j=1,\cdots m$.\\
Moreover for every $\mathscr{B}_b(X)$-normalized eigenfunction $(v_n^{(j)})$ of $\Lambda_n^{(j)}$ there is a subsequence $(v_{n_k}^{(j)})$ which converges uniformly to an eigenvector of $\lam_\infty$.
\end{enumerate}
\label{Eigenv-Convergence}
\end{theo}
The eigenvalues appearing in assertion 2. are called the $\lam$-group of $\lam_\infty$.
\begin{proof}
The first assertion follows from Theorem \ref{Uniform-Convergence}.\\
The proof of the rest follows from assertion 1. and  the general theory of operators approximation in Banach spaces. Indeed, let $\lam_\infty$ be a strictly positive  eigenvalue of $G^{\mu_\infty}$ of multiplicity $m$. By Remark \ref{general-remark-2}, $\lam_\infty$ is also an isolated eigenvalues of $G_{\mathbb C}^{\mu_\infty}$ with the same multiplicity  and hence there is a positively oriented circle  $\Gamma$ enclosing $\lam_\infty$ but no other eigenvalues of $G_{\mathbb C}^{\mu_\infty}$. Upper semi-continuity of the set-valued spectrum function under operator norm convergence implies that for large $n$ the circle $\Gamma$ does not intersect the spectrum of $G_{\mathbb C}^{\mu_n}$. Set
\[
Q^{\mu_\infty} := -\frac{1}{2i\pi} \int_\Gamma (G_{\mathbb C}^{\mu_\infty} - z)^{-1}\,dz,\  
Q^{\mu_n}:= -\frac{1}{2i\pi} \int_\Gamma (G_{\mathbb C}^{\mu_n} - z)^{-1}\,dz.
\]
Then $Q^{\mu_n}$ is well defined for large $n$. Moreover, uniform convergence of $G^{\mu_n}$ towards $G^{\mu_\infty}$ implies uniform convergence of $G_{\mathbb C}^{\mu_n}$ towards $G_{\mathbb C}^{\mu_\infty}$ which in turns  implies that $Q^{\mu_n}\to Q^{\mu_\infty}$ uniformly. It follows in particular that for large $n$, we have (see \cite[Lemma 1.5.5, p.31]{Davies-Linear})
\[
\dim\ran (Q^{\mu_n} )= \dim\ran (Q^{\mu_\infty}) =m.
\]
Thus $\Gamma$ encloses  $\Lambda_n^{(j)}, j=1,\cdots,m$ eigenvalues of $G_{\mathbb C}^{\mu_n}$ and hence of $G^{\mu_n}$, by Remark \ref{general-remark-2}. As $\Gamma$ can be chosen arbitrary small they all converge to $\lam_\infty$.\\
We omit the superscript $j$. Let $v_n$ be a $\mathscr{B}_b(X)$-normalized eigenfunction of $\Lambda_n$. Owing to compactness of $G^{\mu_\infty}$ and $(G^{\mu_\infty} - G^{\mu_n})$ there is a subsequence which we still denote by $(v_n)$, $w\neq 0$ such that $G^{\mu_n} v_n=\Lambda_n v_n\to w$. Thus $v_n\to v:=\frac{1}{\lam_\infty} w\neq 0$ and 
\[
G^{\mu_\infty} v = \lim_{n\to\infty} G^{\mu_\infty} v_n =    \lim_{n\to\infty} (G^{\mu_\infty}  -  G^{\mu_n}) v_n  + \lim_{n\to\infty} G^{\mu_n} v_n = w = \lam_\infty v.  
\] 

\end{proof}
For each $n\in\N_\infty$ let $E_n^{(k)} := \frac{1}{ \lam_n^{(k)} }$ be the enumeration of the eigenvalues of $\calE^{\mu_n}$: 
\[
E_n^{(1)}\leq  E_n^{(2)}\leq\cdots\leq E_n^{(k)}\leq\cdots(\to\infty).
\]
Here is a second aspect of continuous dependence of the trace form with respect to the measure.
\begin{lem}
Let $\lam>0$, and for any interval $I\subset\R$ let $P_n(I), n\in\N_\infty$ be the sepctral projection with respect to $L_{\mu_n}$. Set $H_n:=L^2(X,\mu_n), n\in\N_\infty$. Then for large integers $n$ it holds
\begin{equation}
\dim P_n(-\infty,\lam)H_n = \dim P_\infty(-\infty,\lam)H_\infty.
\end{equation} 
\label{constant-dim}
\end{lem}
\begin{proof}
Let $\lam>0$.\\
{\em Case 1}: $\lam$ is not an eigenvalue of $L_{\mu_\infty}$. Owing to the spectral properties of $L_{\mu_\infty}$, the latter  has only a finite number of strictly positive eigenvalues  with finite multiplicities  in $(-\infty,\lam)$, say $E_1,\cdots,E_N$  with $\lam_i:=1/{E_i}$ being eigenvalues of $G^{\mu_\infty}$. Thus there is $a>0$ such that $a<E_i<\lam, i=1,\cdots,N$. Continuity of finite systems of eigenvalues with respect to uniform convergence (applied to $G_{\mathbb C}^{\mu_n}$, see \cite[p.213]{Kato}) together with Proposition \ref{Resolvent-1} and Remark \ref{general-remark-2}, imply that for  large $n$ we have $\sigma(L_{\mu_n})\cap(-\infty,\lam) = \sigma(L_{\mu_n})\cap(a,\lam)$ and the latter consists of a finite number of eigenvalues. Thereby we obtain
\begin{equation}
P_n(-\infty,\lam)H_n  =  P_n(a,\lam)H_n\ \text{for large}\ n.
\label{equality-spec-projection}	
\end{equation}
To continue the proof  we use Lemma \ref{dimension-equal}. Let $\Gamma$  be a positively oriented circle joining $1/a$ and $1/\lam$. Then for $n$ large enough $\Gamma\subset\rho(G_{\mathbb C}^{\mu_n})\setminus\{0\} = \rho(G{^\mu_n})\setminus\{0\}$. 
In the rest of the proof we suppose that $n$ is large enough. Let 
\begin{align*}
Q^{\mu_n} := -\frac{1}{2i\pi} \int_{\Gamma} (G_{\mathbb C}^{\mu_n} - z)^{-1}\,dz.
\end{align*}
Using Lemma \ref{dimension-equal} (for the measure $\mu_\infty$) together with the equality (\ref{equality-spec-projection}) we obtain 
\begin{equation}
\dim P_\infty(-\infty,\lam)H_\infty = \dim P_\infty(a,\lam)H_\infty = \dim \ran (Q^{\mu_\infty}) .
\label{equality-dim1}
\end{equation}
Moreover using Theorem \ref{Eigenv-Convergence} together with (\ref{equality-dim1}), we obtain for large $n$
\begin{align}
\dim \ran (Q^{\mu_n}) = \dim \ran (Q^{\mu_\infty} ) = \dim P_\infty(a,\lam)H_\infty = \dim P_\infty(-\infty,\lam)H_\infty.
\label{equality-dim2}
\end{align}
Once again, using Lemma \ref{dimension-equal} (for the measures $\mu_n$) together with the equality (\ref{equality-spec-projection}) we achieve 
\begin{equation}
\dim \ran (Q^{\mu_n}) =  \dim P_n(a,\lam)H_n  = \dim P_n(-\infty,\lam)H_n\ \text{for large }\ n.
\label{equality-dim3}
\end{equation}
Finally using (\ref{equality-dim1}), (\ref{equality-dim2}) and (\ref{equality-dim3}) we obtain, for large integers $n$ 
\[
\dim P_\infty(-\infty,\lam)H_\infty = \dim P_n(-\infty,\lam)H_n.
\]
%
%
{\em Case 2:} $\lam$ is an eigenvalue of $L_{\mu_\infty}$.  Then it is an  isolated eigenvalue with finite multiplicity and hence there is $\epsilon>0$ such that $\lam\pm\epsilon$ are not eigenvalues of $L_{\mu_\infty}$ and $\dim\ran P_\infty(\lam-\epsilon,\lam+\epsilon)<\infty$. Owing to  upper semi-continuity of the spectrum together with Proposition \ref{Resolvent-1} we conclude that $\lam\pm\epsilon$.\\
Now it suffices to observe that 
\begin{align*}
	\dim P_n(-\infty,\lam)H_n= \dim P_n(-\infty,\lam-\epsilon)H_n + \dim P_n(\lam - \epsilon,\lam + \epsilon)H_n, n\in\N_\infty,
\end{align*}
and apply  the results obtained in the first case in conjunctions with the arguments used in its proof to get the sought equality, which ends the proof.
\end{proof}
We are in position now to resume continuity aspects of the the trace form with respect to the involved measure.
\begin{theo}
We have
\begin{enumerate}
\item Let $E_\infty$ be an eigenvalue of $\calE^{\mu_\infty}$  of multiplicity $m$. Then for large $n$ there is  $m$  eigenvalues (counted as many times as their respective multiplicities) ${\tilde E}_n^{(1)},\cdots\tilde E_n^{(m)}$ of $\calE^{\mu_n}$ such that  $\lim_{n\to\infty} \tilde E_n^{(j)} = E_\infty, j=1,\cdots m$. Moreover for every $\mathscr{B}_b(X)$-normalized eigenfunction $(v_n^{(j)})$ of $\tilde E_n^{(j)}$ there is a subsequence $(v_{n_k}^{(j)})$ which converges uniformly to an eigenvector of $E_\infty$.
\item 
\begin{equation}
		\lim_{n\to\infty} E_n^{(k)} = E_\infty^{(k)},\ \text{for all}\ k.
\label{k-eigen-convergence}		
\end{equation}
\end{enumerate}
\label{Eigenv-Conv-Form}
\end{theo}
\begin{proof}
The first assertion follows from  Proposition \ref{Resolvent-1} in conjunction with  Theorem \ref{Eigenv-Convergence}.\\	
We proceed to prove the second assertion. Let $P_n(\cdot), H_n, n\in\N_\infty$ be as in Lemma \ref{constant-dim}. By the characterization of ordered eigenvalues (see  \cite[Theorem 12.1, p.266]{Schmudgen}) we have
\begin{equation}
	E_n^{(k)} = \inf_{\dim P_n(-\infty,\lam)H_n\geq k}\,\lam.
\label{Schmudgen}
\end{equation}
Let $\beta$ be such that $\dim P_\infty(-\infty,\beta)H_\infty\geq k$. Then by Lemma \ref{constant-dim} we have, for large $n$, $\dim P_n(-\infty,\beta)H_n\geq k$ and hence, using formula (\ref{Schmudgen}) we derive
\[
\lim_{n\to\infty} E_n^{(k)} = \lim_{n\to\infty} \inf_{ \dim P_n(-\infty,\lam)H_n\geq k  } \lam \leq \beta.
\]
Passing to the infimum on all $\beta$ on the right-hand-side and using using formula (\ref{Schmudgen}) once again,  we get $\lim_{n\to\infty} E_n^{(k)}\leq E_\infty^{(k)}$.\\
Conversely let $\lam$ be such that $\dim P_n{(-\infty,\lam)}H_n\geq k$ for large $n$. By Lemma \ref{constant-dim}, once again, we obtain $\dim P_\infty{(-\infty,\lam)}H_\infty\geq k$.  Using formula (\ref{Schmudgen}), once again we obtain
\[
E_\infty^{(k)} \leq  \lim_{n\to\infty} \inf_{ \dim P_n(-\infty,\lam)H_n\geq k  } \lam = \lim_{n\to\infty} E_n^{(k)},
\]
and we are done.
\end{proof}
\begin{rk}
{\rm
Let us emphasize that Theorem \ref{Eigenv-Conv-Form} implies convergence of the inf-sup involved in Rayleigh-Ritz variational formula for ordered eigenvalues. Precisely, we have
\begin{align*}
\lim_{n\to\infty}\inf_{\dim {\mathscr L} = k}\sup\big\{\calE^{\mu_n}[u], u\in {\mathscr L}\subset \dom(\calE^{\mu_n}), \|u\|_{\mu_n} =1\big\} =\\
\inf_{\dim {\mathscr L} = k}\sup\big\{\calE^{\mu_\infty}[u], u\in {\mathscr L}\subset \dom(\calE^{\mu_\infty}), \|u\|_{\mu_\infty} =1\big\}.
\end{align*}

}
\end{rk}
Let us turn our attention to convergence of resolvents. 
\begin{lem}
Let $\mu\in\mathscr{K}_b(X)$. Then
\begin{enumerate}
\item For every $\alpha\geq 0$ the operator $1+\alpha G^\mu$ is invertible with bounded inverse on $\mathscr{B}_b(X)$. Set
\begin{equation}
	\tilde R_\alpha :=  (1 + \alpha G^\mu)^{-1} G^\mu: \mathscr{B}_b(X)\to \mathscr{B}_b(X),\ \alpha\geq 0.
	\label{Restricted-Res}
\end{equation}
\item $\tilde R_\alpha = R_\alpha|_{\mathscr{B}_b(X)} $.
\end{enumerate}
\label{Resolvent-Bb}
\end{lem}
\begin{proof}
1. Since $G^\mu$ is compact it suffices to prove that $1+\alpha G^\mu$ is injective. Let $u\in \mathscr{B}_b(X)$ be such that $u + \alpha G^\mu u =0$. Then $u + \alpha K^\mu u = 0\ \mu-a.e.$ and by invertibility of $1 + \alpha K^\mu$ we conclude that $u=0\ \mu-a.e.$. But then $u = -\alpha\int_X G(\cdot,y)u(y)\,d\mu(y)=0$. Thus $1+\alpha G^\mu$ is injective.\\
Assertion 2. is obvious.

\end{proof}
\begin{prop}
Let $(\mu_n)_{n\in\N_\infty}$ be as in (\ref{Weak-Monotone}), and
\[
\tilde R_\alpha^{(n)}:= (1 + \alpha G^{\mu_n})^{-1} G^{\mu_n},\ n\in\N_\infty,\ \alpha\geq 0.
\]
Then $\tilde R_\alpha^{(n)}\to \tilde R_\alpha^{(\infty)}$, uniformly.
\end{prop}
\begin{proof}
The result is an immediate consequence of Theorem \ref{Uniform-Convergence}.	
\end{proof}	
\subsection{Quantitative results}
For each $n\in\N_\infty$ let $(E_n^{(k)})_k$ be the eigenvalues of $\calE^{\mu_n}$; ordered in an increasing way and repeated as many times as their respective multiplicities. Set $\nu_n:= |\mu_\infty - \mu_n|$.
\begin{theo}
Let $E_\infty$ be an eigenvalue of $\calE^{\mu_\infty}$ and let $\frac{1}{E_n(j)}$ be any element from the $\lambda$-group of $\lam_\infty:=1/{E_\infty}$. Then  there is a finite constant $c>0$ such that for large integers $n$ it holds
\[
\Big|\frac{1}{E_\infty} - \frac{1}{E_n(j)}\Big| \leq c\|G^{\nu_n}1\|_\infty. 
\]
\label{Eigen-Error}
\end{theo}
\begin{proof}
The proof is mainly taken from \cite{Atkinson,Osborn} we reproduce it in this special context for the convenience of the reader.\\
Let $\lam_n (j)= 1/{E_n(j)}$. Denote by $Q^{\mu_\infty}$ the Riesz projection of $\lam_\infty$ and by $Q^{\mu_n}$ the total Riesz projection of its $\lam$-group.\\
As $\lim_{n\to\infty}\| Q^{\mu_n} - Q^{\mu_\infty}\| = 0$, we get $\dim \ran (Q^{\mu_n}) = \dim\ran(Q^{\mu_\infty})$ for large $n$ and hence the operators
\[
Q_n := Q^{\mu_n}|_{\ran(Q^{\mu_\infty)}}: \ran (Q^{\mu_\infty})\to \ran (Q^{\mu_n}),
\] 
define an isomorphism for large $n$ (see the proof of  \cite [Theorem 2]{Osborn}) . From now on we assume that $n$ is large enough.\\ 
Let us quote that (see the proof of  \cite[Theorem 2]{Osborn} for (\ref{inverses} below))
\begin{align}
Q^{\mu_n} Q_n^{-1} u &= u,\ \forall\,u\in\ran(Q^{\mu_n}),\ \text{and}\  Q_n^{-1} Q^{\mu_n} u = u,\ \forall\,u\in\ran(Q^{\mu_\infty})\nonumber\\ 
&\text{and}\ c:=\sup_n  \|Q_n^{-1} Q^{\mu_n}\|<\infty.
\label{inverses}
\end{align}
Set
\[
\hat{G}^{\mu_n}:= Q_n^{-1} G_{\mathbb C}^{\mu_n}Q_n:\ran(Q^{\mu_\infty}) \to \ran(Q^{\mu_\infty}).
\]
Then $\sigma(\hat{G}^{\mu_n})$ coincides with the $\lambda$-group of $\lam_\infty$ and thereby $\lam_n(j)$ is an eigenvalue of $\hat{G}^{\mu_n}$.\\
Let $u$ be a real-valued  eigenfunction of $\hat{G}^{\mu_n}$ corresponding to $\lam_n(j)$ such that $\|u\|_\infty =1$. Taking into account that $u\in\ran(Q^{\mu_\infty})$, using (\ref{inverses}) and the fact that $Q^{\mu_n}$ commutes with $G_{\mathbb C}^{\mu_n}$ and hence with $G^{\mu_n}$,  an elementary computation leads to
\begin{align}
\big(\lam_\infty - \lam_n(j)\big)u & = \big( G^{\mu_\infty} - \lam_n(j) \big) u = Q_n^{-1}Q^{\mu_n}\big( G^{\mu_\infty}u - \lam_n(j) u\big)\nonumber\\
& =Q_n^{-1}Q^{\mu_n} G^{\mu_\infty}u - \lam_n(j) Q_n^{-1}Q^{\mu_n} u \nonumber\\
&  =  Q_n^{-1}Q^{\mu_n} G^{\mu_\infty}u -  Q_n^{-1}Q^{\mu_n} \hat{G}^{\mu_n} u
= Q_n^{-1}Q^{\mu_n} G^{\mu_\infty}u -  Q_n^{-1}Q^{\mu_n} Q_n^{-1} G^{\mu_n} Q_n u\nonumber\\
&= Q_n^{-1}Q^{\mu_n} G^{\mu_\infty}u -    Q_n^{-1}Q^{\mu_n} Q_n^{-1} G^{\mu_n} Q^{\mu_n} u \nonumber\\ 
&= Q_n^{-1}Q^{\mu_n} G^{\mu_\infty}u -    Q_n^{-1}Q^{\mu_n} Q_n^{-1} P_n G^{\mu_n}  u  \nonumber\\
&= Q_n^{-1}Q^{\mu_n} G^{\mu_\infty}u -    Q_n^{-1} Q^{\mu_n} G^{\mu_n}  u  = Q_n^{-1}Q^{\mu_n} \big(G^{\mu_\infty} - G^{\mu_n}\big)u.
\end{align}
Passing to the norm we get
\[
|\lam_\infty - \lam_n(j)| \leq  \|Q_n^{-1}Q^{\mu_n} \|\|G^{\mu_\infty} - G^{\mu_n}\|\leq c \|G^{\nu_n}1\|_\infty ,
\]
which leads to the sought estimate.
\end{proof}
There are instances where the error estimate from Theorem \ref{Eigen-Error} can be improved.
\begin{lem}
For every $\mu\in\mathscr{K}_b(X)$, the Dirichlet form $\calE^\mu$ is irreducible.
\end{lem}
\begin{proof}
It suffices to prove irreducibility of $K^\mu$. Let $A$ be a Borel subset of $X$ such that $K^\mu 1_A = 1_A K^\mu$. Then $K^\mu 1_A =0$ on $A^c$. However, since $G>0$ we conclude that either $\mu(A)=0$ or $\mu(A^c)=0$. Thus either $A=X,\ \mu-a.e.$ or $A^c = X\,\mu-a.e.$ and $K^\mu$ is irreducible.
\end{proof}
\begin{prop}
Let $E_n^{(0)}, n\in\N_\infty$ be the smallest eigenvalue of $\calE^{\mu_n}, n\in\N_\infty$. Then
\begin{enumerate}
	\item The sequence $(E_n^{(0)})_{n\in\N}$ is monotone  decreasing if $\mu_n\uparrow \mu_\infty$ and monotone  increasing if $\mu_n\downarrow \mu_\infty$.
	\item \[
	E_n^{(0)} -  E_\infty^{(0)}  =   \frac{1}{ \|K^{\mu_n}\| } -  \frac{1}{ \|K^{\mu_\infty}\| },\ n\in\N.
	\]
\end{enumerate}	
\end{prop}
\begin{proof}
Owing to irreducibility property,  for each $n\in\N_\infty$ the smallest eigenvalue of  $\calE^{\mu_n}$, $E_n^{(0)}$ is simple with an eigenfunction that can be chosen positive, say $u_n$. Accordingly, with $\lam_n^{(0)}:=(E_n^{(0)})^{-1} $we obtain
\[
\lam_n^{(0)} u_n = G^{\mu_n} u_n \leq G^{\mu_{n+1}} u_n = K^{\mu_{n+1}} u_n,
\] 
Using compactness, selfadjointness and positivity of the the operators $K^{\mu_n}$, we derive
\[
(E_n^{(0)})^{-1} = \lam_n^{(0)}  \leq 
\sup\{  \frac{  \int_X u K^{\mu_{n+1}} u\,d\mu_{n+1}}{\int_X u^2\,d\mu_{n+1} },\ u\in L^2(X,\mu_{n+1})\setminus\{0\} \}
 = \lam_{n+1}^{(0)} = (E_{n+1}^{(0)})^{-1}.
\]
Thereby the sequence $(E_n^{(0)})_n$ is monotone decreasing.\\
Finally, owing to compactness, selfadjointness and positivity of the the operators $K^{\mu_n}$ once again we obtain
\[
 E_n^{(0)} -  E_\infty^{(0)}  =   \frac{1}{ \|K^{\mu_n}\| } -  \frac{1}{ \|K^{\mu_\infty}\| },\ n\in\N.
\]
\end{proof}
\section{Applications}
\label{Applications}
\subsection{Continuity of stationary solutions for Dirichlet operators with measure-valued coefficients}
We first give an application which involves continuity of the resolvent with respect to the measure.\\
\indent Let $\calE$ be the gradient energy form in $L^2(\R^d,dx), d\geq 3$. It is well known that $\calE$ is transient and has a strictly positive jointly continuous  heat kernel. Moreover it is related to the Laplace operator on  $\R^d$.\\
Let $(\mu_n)_{n\in\N_\infty}\subset\mathscr{K}_b(\R^d)$ be such that $\mu_n\to\mu_\infty$ weakly-monotonically. In these circumstances we have, for all $n$
\begin{equation}
G^{\mu_n} u (x) = c_d\int_{\R^d} \frac{u(y)}{|x - y|^{d-2}}\,d\mu_n(y), u\in\mathscr{B}_b(\R^d), x\in\R^d. 
\label{Green-Pot}
\end{equation}
\begin{prop}
Let $u\in\mathscr{B}_b(\R^d)$ and $\alpha>0$. Then,	
\begin{enumerate}
\item For any $n\in\N_\infty$ the equation 
\begin{equation}
-\Delta u_n + \alpha u_n\mu_n =  u\mu_n\ \text{in the sense of distributions on}\ \R^d,
\label{stationary}
\end{equation}
has a unique solution in $u_n\in C_0(\R^d)$.
\item The solution is given by $u_n = (1 + \alpha G^{\mu_n})^{-1}G^{\mu_n} u$.
\item It holds
\[
\|u_n - u_\infty\|_\infty \leq c_d\|u\|_\infty\cdot\sup_{x\in\R^d}\big( \int_{\R^d} |x - y|^{2-d}\,d(|\mu_\infty - \mu_n|)(y)\big).
\]
\item $\|u_n - u_\infty\|_\infty\to 0$ as $n\to\infty$.
\end{enumerate}
\label{Conv-Stat}		
\end{prop}	
\begin{proof}
Let $u_n,v_n\in C_0(\R^d)$ be two solutions of equation (\ref{stationary}). Then with $w_n:=u_n - v_n$ we have $w_n\in C_0(\R^d)$ is continuous  and $-\Delta w_n + \alpha w_n\mu_n = 0$ in the sense of distributions. Owing to the identity 
\begin{equation}
-\Delta G^\mu f = f\mu\ \text{in the sense of distributions}, \mu\in\mathscr{K}_b(\R^d),f\in\mathscr{B}_b(\R^d),
\label{Potential_Identity}
\end{equation}
 we obtain $-\Delta (w_n + \alpha G^{\mu_n} w_n) = 0$ in the sense of distributions. Elliptic regularity (hypoellipticity of the Laplacian together with Liouville property) implies that $w_n + \alpha G^{\mu_n} w_n$ is constant. Now using Lemma \ref{Kato-property} we get $w_n =0$ which leads to uniqueness of the solution.\\
By Lemma \ref{Kato-property} once again, $u_n:= (1 + \alpha G^{\mu_n})^{-1}G^{\mu_n} u = G^{\mu_n}(1 + \alpha G^{\mu_n})^{-1}u\in C_0(\R^d)$ and satisfies $u_n +\alpha G^{\mu_n} u_n = G^{\mu_n}u$. Applying the Laplacian on both sides and using identity (\ref{Potential_Identity}) we achieve that $u_n$ is the solution.\\
A straightforward computation leads to 
\[
u_\infty - u_n = (1 + \alpha G^{\mu_n})^{-1}(G^{\mu_\infty} - G^{\mu_n})(1 + \alpha G^{\mu_n})^{-1}u.
\]
Using the fact that $\| (1 + \alpha G^{\mu_n})^{-1}\|\leq 1$ in conjunction with Theorem \ref{Uniform-Convergence} and formula (\ref{Green-Pot}) we obtain the sought estimate. Finally assertion 4 follows from 3.
\end{proof}
\begin{rk}
{\rm
\begin{enumerate}
\item One can replace $\R^d$ by a Riemannian manifold $M$ and the Laplace by the Laplace--Beltrami operator on $M$ provided $M$ has a 'nice' heat kernel and still get the same results.	
\item One can replace the Laplacian on $\R^d$ by any elliptic symmetric operator having a 'nice' heat kernel and still get the same results.
\item One can replace the gradient energy form by the Dirichlet form corresponding to the fractional Laplacian and get the same results.	
\end{enumerate}
}
\end{rk}
\subsection{Approximation of the eigenvalues of perturbed one-dimensional graph-Laplacian}

Let $X=\R, m=dx$ the Lebesgue measure on $\R$ and
\[
\calD= H^1(\R), \calE[u] =\int_\R (u'(x))^2\,dx +\int_\R u^2(x)\,dx.
\]
In this situation we have $G(x,y)= \frac{1}{2} e^{-|x-y|}, x,y\in\R$.\\ 
\indent Let $(a_k)_{k\in\Z}\subset (0,\infty)$ be such that $\sum_k a_k <\infty$ and $\mu_n =\sum_{|k|\leq n} a_k\delta_k$. Then $\mu_n\uparrow \mu_\infty$ and by Sobolev inequality each $\mu_n$ satisfies Hardy's inequality. Besides it is easy to check that $(\mu_n)_n\subset\mathscr{K}_b(\R)$. For each $n\in\N_\infty$ set $F_n$ the support of $\mu_n$. Then
\[
J_n :H^1(\R)\to L^2(\R,\mu_),\ u\mapsto u|_{F_n}, K^{\mu_n} u =\frac{1}{2}\int_\R e^{-|\cdot-y|}u(y)\,d\mu_n(y), u\in  L^2(\R,\mu_n), n\in\N_\infty.
\]
Thus the $J_n$'s have dense range and $K^{\mu_n}$'s are compact and hence $\calE^{\mu_n}$ 's have compact resolvent.\\
\indent Let us compute $\calE^{\mu_n}$. A computation similar to the one made in \cite[Example 4.3]{BBST}  leads to:\\
For $n\in\N$ it holds $\dom (\calE^{\mu_n}) = \ran (J_n) = \R^{2n+1}$ and

\begin{align*}
	\calE^{\mu_n}(Ju,Ju) & =\frac{1}{\sinh 1} \sum_{|k|\leq n-1} (u(k+1)-u(k))^2 + 2\frac{\cosh 1 - 1}{\sinh 1} \sum_{|k|\leq n-1} u^2(k)\\
&+ ( u^2(-n) + u^2(n) ).
\end{align*}
Whereas $\dom (\calE^{\mu_\infty})=\ran (J_\infty) = \bigl\{\psi\in L^2(\R,\mu):\; \sum_{n\in\Z}|\psi(n)|^2<\infty\bigr\}$.

\begin{align*}
	\calE^{\mu_\infty}(Ju,Ju) = \frac{1}{\sinh 1} \sum_{k\in\Z} (u(k+1)-u(k))^2 + 2\frac{\cosh 1 - 1}{\sinh 1} \sum_{k\in\Z} u^2(k).
\end{align*}
Let us quote that for finite $n$ the form $\calE^{\mu_n}$ has finitely many eigenvalues, whereas  $\calE^{\mu_\infty}$ has infinitely many distinct eigenvalues.\\
In this case applying Theorem \ref{Eigenv-Conv-Form} in conjunction with Theorem \ref{Eigen-Error} we obtain
\[
{E_n^{(k)}}\to {E_\infty^{(k)}}\ \text{and}\  |\frac{1}{E_\infty^{(k)}} - \frac{1}{E_n^{(k)}}| \leq c\sup_{x\in\R}\sum_{|j|>n} a_j e^{-|x - j|}.
\]

\subsection{From the final eigenvalues to the approximating eigenvalues}
Here we give an example where eigenvalues of the final form are known and those of the approximating trace forms are not explicitly known. It is about substituting eigenvalues of D-to-N operators  on 'small' sets by those corresponding to D-to-N operators on thin sets, namely boundaries.\\
\indent Let $\Om$ be an open non-empty connected bounded subset of $\R^d,\ d\geq 3$ having the extension property (for example smooth domain), $\Gamma:=\partial\Om$ and $\Om_n:=\{ x\in\Om\colon dist(x,\Gamma)<1/n\}, n\in\N$. Take $X=\overline\Om$ and let us  consider the regular transient Dirichlet form  $\calE$:
\[
\dom(\calE) = H^1(\Om)\subset L^2(\overline \Om,dx),\ \calE[u] = \int_\Om |\nabla u|^2\,dx + \int_\Om u^2\,dx.
\]
It is well known that the semigroup associated with $\calE$ has a symmetric  kernel $p_t(\cdot,\cdot)$ which is jointly continuous, non-negative (see \cite{Davies}). Besides, according to \cite[Theorem 6.10, p.171]{Ouhabaz}, $p_t$ satisfies the upper Gaussian estimate
\begin{equation}
p_t(x,y) \leq \frac{c}{t^{d/2}} e^{- \frac{|x - y|^2}{4t}},\ t>0,x,y\in X.	
\end{equation}
These results are mainly derived from uniform ellipticity and Sobolev embeddings.\\	 
Continuity of $p_t$ in conjunction with formula (\ref{laplace-transform}) and Fatou Lemma yield that $G(\cdot,y)$ is l.s.c. for any $y\in X$.\\
Continuity of $p_t$, once again together with the upper Gaussian estimate imply  that $G$ is continuous away from the diagonal and 
\begin{equation}
G(x,y) \leq c|x-y|^{2 - d},\ x,y\in X.
\label{Green-bound}
\end{equation}
\indent Let $\mu_\infty$ be the surface  measure of  $\Gamma$, i.e. the $(d-1)$-Hausdorff measure of $\Gamma$. By \cite[Example 3, p.30]{Jonsson-Wallin} we have
\begin{equation}
\mu_\infty(B_x(r)) \sim r^{d-1},\ x\in\Gamma, 0<r\leq 1. 
\label{Vd-Ball}
\end{equation}
Hence,  (\ref{Green-bound}), (\ref{Vd-Ball}) and Corollary \ref{Vol-Growth} imply that $\mu_\infty$ is $G$-Kato and since it is finite we obtain $\mu_\infty\in\mathscr{K}_b(X)$.
Let  
\[
\mu_n = 1_{\Om_n}dx + \mu_\infty, n\in\N.
\]
One can check that $\mu_n\downarrow dS$ weakly and arguing as before one get  $(\mu_n)_n\subset \mathscr{K}_b(X)$. Moreover for each $n$ the support $F_n$ of $\mu_n$ is $\overline\Om_n$.\\
\indent We choose $J_n, n\in\N$ as follows:
\[
\dom (J_n) = H^1(\Om)\to L^2(X,\mu_n), u\mapsto u|_{F_n} =  u|_{\Om_n} + tr_\Gamma u,
\]   
and
\[
\dom (J_\infty) = H^1(\Om)\to L^2(X,\mu_\infty), u\mapsto u|_{\Gamma} =  tr_\Gamma u,
\] 
where $tr_\Gamma u$ is the trace to the boundary $\Gamma$ of $u$.\\
Then the form $\calE^{\mu_\infty}$ is given by
\begin{eqnarray*}
	\dom(\calE^{\mu_\infty})= H^{1/2}(\Gamma),\ \calE^{\mu_\infty}[\psi] = \calE[Pu]=\int_\Om |\nabla P u|^2\,dx + \int_\Om (P u)^2\,dx\,\ \forall\,\psi\in H^{1/2}(\Gamma),
\end{eqnarray*}
where $Pu$ is the unique solution in $H^1(\Om)$ of 
\begin{eqnarray*}
	\left\{\begin{gathered}
		-\Delta P u + P u=0, \quad \hbox{in } \Om,\\
		P u= \psi,~~~{\rm on}\  \Gamma
	\end{gathered}
	\right.
\end{eqnarray*}
Using Green's formula (for further details we refer to \cite[section 6.2]{BenamorJMAA} ) we achieve
\[
\calE^{\mu_\infty}[J_\infty u] = \int_\Gamma \frac{\partial Pu}{\partial\nu} u\,d\mu_\infty,\ u\in H^1(\Om),
\]
where $\frac{\partial \cdot}{\partial\nu}$ is the outer normal derivative.\\
\indent For $n\in\N$ we denote by $\Gamma_n$ the part of the boundary of $\Om_n$ which lies inside $\Om$ and $dS_n$ its surface measure. Similar computations lead to
\[
\calE^{\mu_n}[J_n u] = \int_{\Om_n} |\nabla u|^2\,dx +   \int_{\Om_n} u^2\,dx  +
  \int_\Gamma \frac{\partial P_n u}{\partial\nu_n} u\,d\mu_\infty
- \int_{\Gamma_n} \frac{\partial P_nu}{\partial\nu} u\,dS_n ,\ u\in H^1(\Om),
\]
where $P_n u$ is the unique solution in $H^1(\Om_n)$ of 
\begin{eqnarray*}
	\left\{\begin{gathered}
		-\Delta P_n u + P_n u=0, \quad \hbox{in } \Om\setminus\overline \Om_n,\\
		P_n u= u,~~~{\rm on}\  \overline \Om_n
	\end{gathered}
	\right.
\end{eqnarray*}
Also for each $n\in\N_\infty$ we have
\begin{align*}
&G^{\mu_n} u = \int_X G(x,y)u(y)\,d\mu_n(y),\ u\in\mathscr{B}_b(X)\\
& \text{and}\ K^{\mu_n} u = \int_X G(x,y)u(y)\,d\mu_n(y),\ u\in L^2(X,\mu_n).
\end{align*}
\begin{rk}
{\rm
A straightforward computation shows that the positive  self-adjoint operator associated with $\calE^{\mu_n}$ is the shifted Laplacian $-\Delta + 1$ on $\Om_n$ with mixed Neumann-Wentzell boundary conditions:
\[
\frac{\partial u}{\partial\nu_n} = 0\ \text{on}\ \Gamma_n\ \text{and}\ \frac{\partial u}{\partial\nu}  - \Delta u= 0\ \text{on}\ \Gamma.
\]	

}
\end{rk}
\begin{prop}
Let $E_n^{(k)}$ be the $k$-th eigenvalue of $\calE^{\mu_n}, n\in\N_\infty$. Then 
\[
	\lim_{n\to\infty} E_n^{(k)} = E_\infty^{(k)}\ \text{and}\  |\frac{1}{E_n^{(k)}} - \frac{1}{E_\infty^{(k)}} | \leq c/n. 
\]	
\end{prop}
\begin{proof}
Convergence of eigenvalues follows from Theorem \ref{Eigenv-Conv-Form}.\\
Let us  prove the rest. Using the upper estimate for the Green function (\ref{Green-bound}) an elementary computation leads to
\begin{align*}
\|G^{\mu_n} 1 - G^{\mu_\infty} 1\|_\infty = \sup_{x\in X}\int_{\Om_n} G(x,y)\,dx 
\leq c \sup_{x\in X}\int_{\Om_n} |x-y|^{2-d}\,dx 
\leq c/n.
\end{align*}
Now the result follows from Theorem \ref{Eigen-Error}.
\end{proof}
For the particular case $\Om = B$, the open unit ball in $\R^3$, according to \cite[Formula 6.20]{BenamorJMAA}, the eigenvalues of $\calE^{\mu_\infty}$ are given by
\begin{eqnarray}
\big\{ m + 2\sum_{k= 0}^\infty\frac{1}{1 + j_{mk}^2}\colon m\in\N\cup\{0\} \big\},
\label{MLBall}
\end{eqnarray}
where $j_{mk}^2$ are the eigenvalues of the Dirichlet-Laplacian on the unit ball. Each eigenvalue has multiplicty $2m + 1$.\\
Hence applying the latter proposition we obtain:
\begin{prop}
Let $E_n^{(k)}$ be the $k$-th eigenvalue of $\calE^{\mu_n}, n\in\N$. Then there is $m=m(k)$ such that
\[
\lim_{n\to\infty} E_n^{(k)} = E_\infty^{(k)}:= m + 2\sum_{l= 0}^\infty\frac{1}{1 + j_{ml}^2}.
\]	
\end{prop}
\section{Appendix: Tests for $G$-Kato property}
Let $\rho$ be the metric of $X$.
\begin{theo}
Assume that $G$ is continuous on $\{x\neq y\}$ and there is $\beta>0, r_0>0$ and $c>0$ such that
\[
G(x,y) \leq c \rho^{-\beta}(x,y)\ \text{for}\ \rho(x,y)<r_0.
\]
Let $\mu$ be a finite Radon measure on $X$. Then $\mu\in\mathscr{K}_b(X)$ if and only if
\[
\lim_{r\to 0} \sup_{x\in X} \int_{B_r(x)} \rho^{-\beta}(x,y)\,d\mu(y) = 0.
\]
\label{Metric-condition}
\end{theo}
\begin{proof}
Let
\[
g_n(x) = \int_{\{\rho(x,y)>1/n\}} G(x,y)\,d\mu(y), x\in X, n\in\N
\]
By assumptions $(g_n)\subset C_0(X)$ for large $n$. Moreover
\[
G^\mu 1(x) - g_n(x) = \int_{\{\rho(x,y)\leq 1/n\}} G(x,y)\,d\mu(y)\leq c \int_{\{\rho(x,y)\leq 1/n\}} \rho^{-\beta}(x,y)\,d\mu(y)
\] 
hence
\[
\sup_{x\in X}|G^\mu 1(x) - g_n(x)| \leq c \sup_{x\in X}\int_{\{\rho(x,y)\leq 1/n\}} \rho^{-\beta}(x,y)\,d\mu(y)\to 0\ \text{as}\ n\to\infty,
\] 
which completes the proof.
\end{proof}
\begin{coro}
Under assumption of Theorem \ref{Metric-condition}, assume further that there is $s>\beta$ and $c'$ such that
\[
\mu(B_r(x)) \leq c'r^s, x\in X, 0<r\leq r_0.
\]
Then if $\mu$ is finite, it is in $\mathscr{K}_b(X)$.
\label{Vol-Growth} 
\end{coro}
\begin{proof}
The result follows from the identity
\[
\int_{B_\epsilon(x)} \rho^{-\beta}(x,y)\,d\mu(y) = \beta\int_0^\epsilon \frac{\mu(B_r(x))}{r^{\beta}}\,\frac{dr}{r}
+ \frac{\mu(B_\epsilon(x))}{\epsilon^{\beta}}.
\]
\end{proof}
\begin{exa}
{\rm
Let $X=\R^d$ and  $0<\alpha<\min(2,d)$. Let $\calE^{(\alpha)}$ be the Dirichlet form associated with the fractional Laplacian in $L^2(\R^d,dx)$:
\[
\dom(\calE^{(\alpha)}) = W^{\alpha/2,2}(\R^d),\ \calE^{(\alpha)}[u] = c_{d,\alpha}\int_{\R^d} |\xi|^\alpha |\hat u(\xi)|^2\,d\xi.
\]	
Then $\calE^{(\alpha)}$	is transient and has a Green kernel
\[
G^{(\alpha)}(x,y) = \kappa_{d,\alpha} |x - y|^{d-\alpha}.
\]
Let $\mu$ be a positive finite Radon measure on $\R^d$ such that
\[
 \lim_{r\to 0}\sup_{x\in\R^d}\int_{\{|x - y|<r\}}|x - y|^{\alpha-d} \,d\mu(y) =0.
\]
Owing to Theorem \ref{Metric-condition} we get $\mu\in\mathscr{K}_b(\R^d)$.\\
It follows in particular that any finite Radon measure such that 
\[
\mu(B_r(x)) \leq c'r^s, x\in X, 0<r\leq r_0,
\]
with $d-\alpha<s\leq d$ is in $\mathscr{K}_b(\R^d)$.
}
\end{exa}

\bibliography{Biblio-Cont-Trace}

\end{document}